\documentclass[preprint,12pt]{elsarticle}

\usepackage{amsmath,amstext,amssymb,amsfonts}
 \usepackage{amsthm}

\usepackage{graphicx}

\usepackage{graphics,color}
\newtheorem{lemma}{Lemma}[section]
\newtheorem{theorem}[lemma]{Theorem}

\newtheorem{definition}[lemma]{Definition}

\newtheorem{example}[lemma]{Example}

\newtheorem{remark}[lemma]{Remark}

\newcommand{\ilim} {\mathop{\rm lim\,inf}}

\newcommand{\Xx}{\mathtt{X}}

\newcommand{\Aa}{\mathtt{A}}
\newcommand{\Bb}{\mathtt{B}}
\newcommand{\ff}{\mathtt{f}}
\newcommand{\vv}{\mathtt{v}}
\newcommand{\Phii}{\mathtt{\Phi}}

\newcommand{\K}{\mathbb{K}}

\newcommand{\Y}{\mathbb{Y}}

\newcommand{\X}{\mathbb{X}}

\newcommand{\R}{\overline{\mathbb{R}}}

\newcommand{\Gr}{{\rm Gr}}

\journal{}

\begin{document}

\begin{frontmatter}



\title{An example showing that $A$-lower semi-continuity is essential for minimax continuity theorems
}


\author{Eugene~A.~Feinberg}

\address{Department of Applied Mathematics and
Statistics,\\
 Stony Brook University,\\
Stony Brook, NY 11794-3600, USA\\  {eugene.feinberg@sunysb.edu}}

\author{Pavlo~O.~Kasyanov}

\address{Institute for Applied System Analysis,\\
National Technical University of Ukraine\\ ``Igor Sikorsky Kyiv Polytechnic
Institute'',\\ Peremogy ave., 37, build, 35, 03056, Kyiv, Ukraine \\ {kasyanov@i.ua}}

\author{Michael~Z.~Zgurovsky}

\address{National Technical University of Ukraine \\``Igor Sikorsky Kyiv Polytechnic
Institute'',\\ Peremogy ave., 37, build, 1,\\ 03056, Kyiv, Ukraine \\ {zgurovsm@hotmail.com}}

\begin{abstract}
Recently Feinberg et al.~\cite{FKZ2017} established results on continuity properties of minimax values and solution sets for a function of two variables depending on a parameter.  Such minimax problems appear in games with perfect information, when the second player knows the move of the first one, in turn-based games, and in robust optimization. Some of the results in \cite{FKZ2017} are proved under the assumption that the multifunction, defining the domains of the second variable, is $A$-lower semi-continuous.  As shown in \cite{FKZ2017}, the $A$-lower semi-continuity property  is stronger than lower semi-continuity, but in several important cases these properties coincide.  This note provides an example demonstrating that in general the $A$-lower semi-continuity assumption cannot be relaxed to lower semi-continuity.
\end{abstract}


\end{frontmatter}

\section{Introduction}\label{s1}
Recently Feinberg et al.~\cite{FKZ2017} established results on continuity properties of minimax values and solution sets for
a function of two variables depending on a parameter, when decision sets may not be compact. Such minimax values appear in games with perfect information, when the second player knows the move of the first one, in turn-based games, and in robust optimization. Some of the results in \cite{FKZ2017} hold under the assumption that a multifunction defining decision sets of the second player is $A$-lower semi-continuous.  The $A$-lower semi-continuity property of a multifunction was introduced in \cite{FKZ2017}, and it is stronger than lower semi-continuity. However, as shown in \cite{FKZ2017}, these two conditions are equivalent in the following two important cases: (i) decision sets for the second player
do not depend on the first variable, as this takes place in games with simultaneous moves, and (ii) the multifinction defining  decision sets of the first player is upper semi-continuous and compact-valued.  This note provides an example when the corresponding continuity properties of minimax fail when the $A$-lower semi-continuity assumption is relaxed to lower semi-continuity.

Let $\R:=\mathbb{R}\cup\{\pm\infty\}$ and $\mathbb{S}$ be a metric space.
For a nonempty set $S\subset \mathbb{S},$  the notation $f:S\subset\mathbb{S}\mapsto {\R}$ means that for each $s\in S$ the value $f(s)\in{\R}$ is defined.
In general, the function $f$ may be also defined outside of $S.$   The notation $f:\mathbb{S}\mapsto {\R}$ means that the function $f$ is defined on the
entire space $\mathbb{S}.$
This notation is equivalent to the notation $f:\mathbb{S}\subset\mathbb{S}\mapsto {\R}, $ which we do not write explicitly.
For a function $f:S\subset\mathbb{S}\mapsto {\R}$ we sometimes consider its restriction $f\big|_{\tilde{S}}:\tilde{S}\subset\mathbb{S}\mapsto {\R}$ to the set $\tilde{S}\subset S.$
Throughout the note we denote by $\K(\mathbb{S})$ the
\emph{family of all nonempty compact subsets} of ${\mathbb{S}}$
and by $S(\mathbb{S})$ the \emph{family of all nonempty
subsets} of~$\mathbb{S}.$

We recall that, for a
nonempty set $S\subset \mathbb{S},$ a function $f:S\subset \mathbb{S}\mapsto\R$ is
called \textit{lower semi-continuous at $s\in S$}, if for each sequence
$\{s_n\}_{n=1,2,\ldots}\subset S,$ that converges to $s$ in $\mathbb{S},$ the
inequality $\ilim_{n\to\infty} f(s_n)\ge f(s)$ holds. A function
$f:S\subset \mathbb{S}\mapsto\R$ is called \textit{upper semi-continuous at
$s\in S$}, if $-f$ is lower semi-continuous at $s\in S.$ A function
$f:S\subset \mathbb{S}\mapsto\R$ is called \textit{lower / upper semi-continuous} if $f$ is lower / upper semi-continuous at each $s\in S.$ A function
$f:S\subset \mathbb{S}\mapsto\R$ is called \textit{inf-compact on
$S$}, if all the level sets $\{s\in S \, : \,  f(s)\le \lambda\},$ $\lambda\in\mathbb{R},$ are compact in $\mathbb{S}.$
A function
$f:S\subset \mathbb{S}\mapsto\R$ is called \textit{sup-compact on
$S$}, if $-f$ is inf-compact on~$S.$

Let $\X$ and $\Y$ be nonempty sets. For a multifunction $\Phi:\X\mapsto 2^{\Y},$ let
$ {\rm Dom\,}\Phi:=\{x\in\X\,:\, \Phi(x)\ne \emptyset\}.$
A multifunction $\Phi:\X\mapsto 2^{\Y}$ is called \textit{strict} if ${\rm Dom\,}\Phi=\X,$ that is,  $\Phi:\X\mapsto S(\Y)$ or, equivalently,
$\Phi(x)\ne \emptyset$ for each $x\in \X.$
For $Z \subset \X$ define the \textit{graph} of a multifunction
$\Phi:\X\mapsto 2^{\Y},$ restricted to~$Z$:
\[
{\rm Gr}_Z(\Phi)=\{(x,y)\in Z\times\Y\,:\, x\in{\rm Dom\,}\Phi,\, y\in \Phi(x)\}.
\]
When $Z=\X,$ we use the standard notation $\Gr(\Phi)$ for the graph of $\Phi:\X\mapsto 2^{\Y}$ instead of $\Gr_{\X}(\Phi).$

\section{Basic definitions and facts}\label{s2}
Let $\Xx, \Aa$ and $\Bb$ be metric spaces, $\Phii_{\Aa}:\Xx\mapsto S(\Aa)$ and $\Phii_{\Bb}:\Gr(\Phii_{\Aa})\subset\Xx\times\Aa\mapsto S(\Bb)$ be multifunctions, and $\ff:{\rm Gr}(\Phii_{\Bb})\subset \Xx \times \Aa\times\Bb \mapsto \overline{\mathbb{R}}$ be a function. Define
the \textit{worst-loss function}
\begin{equation}\label{eq1starworstloss}
\ff^\sharp(x,a):=\sup\limits_{b\in \Phii_{\Bb}(x,a)}\ff(x,a,b),\qquad
(x,a)\in\Gr(\Phii_{\Aa}),
\end{equation}
the \textit{minimax value function}
\begin{equation}\label{eq1starminimax}
\vv^\sharp(x):=\inf\limits_{a\in \Phii_{\Aa}(x)}\sup\limits_{b\in \Phii_{\Bb}(x,a)}\ff(x,a,b),\qquad
x\in\Xx,
\end{equation}
and the \textit{solution multifunctions}
\begin{equation}\label{e:defFi*minimax1}
\Phii_{\Aa}^*(x):=\big\{a\in
\Phii_{\Aa}(x)\,:\,\vv^\sharp(x)=\ff^\sharp(x,a)\big\},\quad x\in\Xx;
\end{equation}
\begin{equation}\label{e:defFi*minimax2}
\Phii_{\Bb}^*(x,a):=\big\{b\in
\Phii_{\Bb}(x,a)\,:\,\ff^\sharp(x,a)=\ff(x,a,b)\big\},\ (x,a)\in\Gr(\Phii_{\Aa}).
\end{equation}

Formulae (\ref{eq1starworstloss}--\ref{e:defFi*minimax2}) describe the value functions and solution multifunctions for one-step zero-sum games with perfect information, where $\Xx$ is the state space, $\Aa$ and $\Bb$ are the action sets of Players I and II respectively.  Player I knows the state $x$ and selects an action from the set $\Phii_{\Aa}(x).$ Player II knows the state $x$ and the move $a$ chosen by Player I and selects an action $b$ from the set $\Phii_{\Bb}(x,a).$ Then Player I pays Player II the amount $\ff(x,a,b).$  Turn-based games can be usually reduced to games with perfect information.  These formulae also describe a model for robust optimization.  In this case the goal is to choose an action $a$ to minimize possible losses $\ff(x,a,b)$ under the worst possible outcome of the uncertain parameter $b.$

Natural continuity properties of  functions~(\ref{eq1starworstloss},\ref{eq1starminimax}) and solution multifunctions (\ref{e:defFi*minimax1},\ref{e:defFi*minimax2}) are described in \cite{FKZ2017}. The results in \cite{FKZ2017} generalize Berge's theorem and Berge's maximum theorem for possibly noncompact action sets from \cite{Feinberg et al} and \cite{FKV} to minimax settings.  We start with the descriptions of these theorems.

Let $\X$ and $\Y$ be metric spaces. A set-valued mapping ${F}:\X
\mapsto 2^{\Y}$ is \textit{upper semi-continuous} at $x\in{\rm Dom\,}F$ if, for
each neighborhood $\mathcal{G}$ of the set $F(x),$ there is a
neighborhood of $x,$ say $U(x),$ such that
$F(x^*)\subset \mathcal{G}$ for all $x^*\in U(x)\cap {\rm Dom\,}F;$ a
set-valued mapping ${F}:\X
\mapsto 2^{\Y}$ is \textit{lower
semi-continuous} at $x\in{\rm Dom\,}F$ if, for each open set $\mathcal{G}$
with $F(x) \cap \mathcal{G} \neq \emptyset,$ there is a
neighborhood of $x,$ say $U(x),$ such that if $x^*\in
U(x)\cap {\rm Dom\,}F ,$ then
$F(x^*)\cap \mathcal{G}\ne\emptyset$ (see e.g., Berge
\cite[p.~109]{Ber} or Zgurovsky et al. \cite[Chapter~1,
p.~7]{ZMK1}). A set-valued mapping is called \textit{upper / lower
semi-continuous}, if it is upper /
lower
semi-continuous at all $x\in{\rm Dom\,}F.$

Let $\Phi:\X\mapsto 2^\Y$  be a multifunction with ${\rm Dom\,}\Phi\ne\emptyset$, and $u:{\rm Gr}(\Phi)\subset \X \times \Y \mapsto \R$ be a function. Define
the
 \textit{value function}
\[
v(x):=\inf\limits_{y\in \Phi(x)} u(x,y),\qquad
x\in\X,
\]
and the \textit{solution multifunctions}
\[
\Phi^*(x):=\big\{y\in
\Phi(x)\,:\,v(x)=u(x,y)\big\},\quad x\in\X.
\]

First, we formulate two classic facts.
\begin{theorem}\label{thB} {\rm (Berge's theorem; Berge~\cite[Theorem~2, p.~116]{Ber},
		Hu and
Papageorgiou~\cite[Proposition~3.3, p.~83]{Hu})}. If
$u:\X\times\Y\to\overline{\mathbb{R}}$ is a lower semi-continuous
function and $\Phi:\X\to \K(\Y)$ is an upper semi-continuous multifunction, then the
function $v:\X\to\overline{\mathbb{R}}$ is lower semi-continuous and the solution sets $\Phi^*(x)$ are nonempty and compact for all $x\in\X.$
	\end{theorem}
\begin{theorem}\label{thBM} {\rm (Berge's maximum theorem; Berge~\cite[p. 116]{Ber}, Hu and
Papageorgiou~\cite[Theorem~3.4, p.~84]{Hu})}
		If $u:\X\times \Y \mapsto \mathbb{R}$ is a continuous function and  $\Phi:\,\X\mapsto \K(\Y)$ is a continuous multifunction, then the value function $v:\X\to\mathbb{R}$ is continuous and the solution multifunction $\Phi^*:\,\X\mapsto \K(\Y)$ is upper semi-continuous.
	\end{theorem}

Second, we formulate Berge's theorem and Berge's maximum theorem for possibly noncompact sets $\Phi(x).$
\begin{definition} \rm{(Feinberg et al.~\cite[Definition 1.1]{Feinberg et al}, \cite[Definition~1]{FKZ2017}).}\label{def:Kinfcom}
A function $u:\Gr(\Phi)\subset\X\times \Y\mapsto \overline{\mathbb{R}}$ is called
\textit{$\K$-inf-compact} on ${\rm Gr}(\Phi),$
if for every $C\in \K({\rm Dom\,}\Phi)$ this function is inf-compact on ${\rm
Gr}_C(\Phi).$
\end{definition}
In particular, according to \cite[Lemma~3]{FKZ2017}, a function  $u:\Gr(\Phi)\subset\X\times \Y\mapsto \overline{\mathbb{R}}$ is
$\K$-inf-compact on ${\rm Gr}(\Phi)$ in the following two cases: (i) $u:\Gr(\Phi)\subset\X\times \Y\mapsto \overline{\mathbb{R}}$ is an inf-compact function; (ii) the assumptions of Berge's theorem (see Theorem~\ref{thB} above) hold. Note that a function $f:\Gr(\Phi)\subset\X\times \Y\mapsto \overline{\mathbb{R}}$ is called
 \textit{$\K$-sup-compact} on ${\rm Gr}(\Phi)$ if
 the function $-f$ is $\K$-inf-compact on ${\rm
Gr}(\Phi).$ The following lemma provides necessary and sufficient conditions for a function to be $\K$-inf-compact.

\begin{lemma}{\rm(Feinberg et al.~\cite[Lemma~2]{FKZ2017} and Feinberg and Kasyanov \cite[Lemma~2]{FKSVAN})}\label{k-inf-compact-strict}
The function $u:\Gr(\Phi)\subset\X\times \Y\mapsto \overline{\mathbb{R}}$ is
$\K$-inf-compact on ${\rm Gr}(\Phi)$ if and only if the following two assumptions hold:
\begin{itemize}
\item[{\rm
(i)}] $u:\Gr(\Phi)\subset\X\times \Y\mapsto \overline{\mathbb{R}}$ is lower semi-continuous;
\item[{\rm(ii)}] if a sequence $\{x_n \}_{n=1,2,\ldots}$ with values in ${\rm Dom\,}\Phi$
converges in $\X$ and its limit $x$ belongs to ${\rm Dom\,}\Phi,$ then each sequence $\{y_n
\}_{n=1,2,\ldots}$ with $y_n\in \Phi(x_n),$ $n=1,2,\ldots,$ satisfying
the condition that the sequence\\ $\{u(x_n,y_n) \}_{n=1,2,\ldots}$ is
bounded above, has a limit point $y\in \Phi(x).$
\end{itemize}
\end{lemma}
\begin{theorem}\label{MT2} {\rm (Berge's theorem for possibly noncompact decision sets; Feinberg et al.~\cite[Theorem~1]{FKZ2017})}.
If a function $u:{\rm Gr}(\Phi)\subset \X \times \Y \mapsto \overline{\mathbb{R}}$ is
$\K$-inf-compact on $\Gr(\Phi),$ then the value function $v:{\rm Dom\,}\Phi \subset
\X\mapsto \overline{\mathbb{R}}$ is lower semi-continuous.
In addition, the following two properties hold for the solution multifunction  $\Phi^*:$
{\rm(a)} the graph ${\rm Gr}({\Phi}^*)$ is a Borel subset of\, $\X\times \Y;$
{\rm(b)} if $v(x)=+\infty,$
then ${\Phi}^*(x)={\Phi}(x),$ and, if $v(x)<+\infty,$ then
${\Phi}^*(x)\in\K(\Y);$ $x\in {\rm Dom\,}\Phi.$
\end{theorem}

\begin{theorem}\label{MT3} {\rm (Berge's maximum theorem for possibly noncompact decision sets; Feinberg et al.~\cite[Theorems 1.2 and 3.1]{FKV})}.
If   $u: \rm{Gr}(\Phi)\subset \X
\times \Y \to \mathbb{R}$ is a $\K$-inf-compact, upper semi-continuous function on
$\rm{Gr}(\Phi)$ and $\Phi: \X \to S(\Y)$ is a lower semi-continuous multifunction, then   the value function $v: \X\to \mathbb{R}$ is
continuous and the solution multifunction $\Phi^*:\X\to  \K(\Y)$ is 
upper semi-continuous.
\end{theorem}

According to \cite[Lemma~3]{FKZ2017} described above, Theorem~\ref{MT2} is a generalization of Berge's  theorem (Theorem~\ref{thB}), and Theorem~\ref{MT3} is a generalization of Berge's maximum theorem (Theorem~\ref{thBM}). In particular, Theorem~\ref{MT2} is important for inventory control and Markov decision processes; see Feinberg~\cite{Ftut} for details.  Before the notion of $\K$-inf-compactness was introduced in \cite{Feinberg et al},
Luque-V\'asquez and Hern\'andez-Lerma~\cite{LVHL} provided an example of a continuous multifunction $\Phi(x)=\Y$ for all $x\in\X,$  continuous function $u:\X\times\Y\mapsto\mathbb{R},$ such that the function $u(x,\cdot):\Y\to\mathbb{R}$ is inf-compact for all $x\in\X,$ for which the value function $v:\X\mapsto\mathbb{R}$ is not lower semi-continuous.  In this example, the function $u:\X\times \Y\mapsto \mathbb{R}$ is not $\K$-inf-compact on $\X\times\Y.$  This example is used to construct Example~\ref{exa:Alsc} below.

Third, we describe the results on continuity properties of minimax values and solution multifunctions from Feinberg et al.~\cite{FKZ2017}.
We start with the properties that do not use $\Aa$-lower semi-continuity of $\Phii_\Bb;$ see statements~(A,B,C) below.

\begin{definition}\label{defi:uniformAlsc}
{\rm A multifunction $\Phii_{\Bb}:\Gr(\Phii_{\Aa})\subset\Xx\times\Aa\mapsto S(\Bb)$ is called {\it 
$\Aa$-lower semi-continuous}, if the following condition holds:
\begin{itemize}
\item[]
if a sequence $\{x_n\}_{n=1,2,\ldots}$ with values in $\Xx$ converges and its limit $x$ belongs to $\Xx,$ $a_n\in \Phii_\Aa(x_n)$ for each $n=1,2,\ldots,$ and $b\in \Phii_\Bb(x,a)$ for some $a\in \Phii_\Aa(x),$ then there is a sequence $\{b_n\}_{n=1,2,\ldots},$ with $b_n\in \Phii_\Bb (x_n,a_n)$ for each $n=1,2,\ldots,$ such that $b$ is a limit point of the sequence $\{b_n\}_{n=1,2,\ldots}.$
\end{itemize}
}\end{definition}

We recall that  a multifunction $\Phii_{\Bb}:\Gr(\Phii_{\Aa})\subset\Xx\times\Aa\mapsto S(\Bb)$ is 
lower semi-continuous, iff for each $x\in\Xx$ and $a\in\Phii_{\Aa}(x)$, for every sequence $(x_n,a_n)\to (x,a)$ with $x_n\in\Xx   ,$ $a_n\in \Phii_{\Aa}(x_n),$ $n=1,2,\ldots,$ and for every $b\in \Phii_\Bb(x,a),$    there exists a sequence $b_{n}\in \Phii_{\Bb}(x_{n},a_{n})$  such that $b_{n}\to b.$

As follows from the definitions, an $\Aa$-lower semi-continuous multifunction $\Phii_\Bb$ is lower semi-continuous, but the opposite statement is not correct; see Feinberg et al. \cite[Example 5]{FKZ2017}. The following lemma describes two conditions under which  a lower semi-continuous multifunction $\Phii_\Bb$ is $\Aa$-lower semi-continuous.  Case (a) takes place when the first player has compact action sets, and case (b) takes place when
 decision sets for the second player do not depend on the first variable, as this takes place in games with simultaneous moves;
see Ja\'skiewicz and Nowak~\cite{Jan,JanS} and references therein for the results on stochastic games satisfying these conditions.
\begin{lemma}\label{lem:unifAlsc}
Let $\Phii_{\Bb}:\Gr(\Phii_{\Aa})\subset\Xx\times\Aa\mapsto S(\Bb)$ be a lower semi-continuous multifunction. Then the following statements hold:
\begin{itemize}
\item[{\rm(a)}] if $\Phii_{\Aa}:\Xx\mapsto S(\Aa)$ is upper semi-continuous and compact-valued at each $x\in\Xx,$ then $\Phii_{\Bb}:\Gr(\Phii_{\Aa})\subset\Xx\times\Aa\mapsto S(\Bb)$ is  $\Aa$-lower semi-continuous;
\item[{\rm(b)}] if $\Phii_{\Bb}(x,a)$ does not depend on $a\in \Phii_\Aa(x)$ for each $x\in\Xx,$ that is, $\Phii_{\Bb}(x,a_*)=\Phii_{\Bb}(x,a^*)$ for each $(x,a_*),(x,a^*)\in \Gr(\Phii_\Aa),$ then $\Phii_{\Bb}:\Gr(\Phii_{\Aa})$ $\subset \Xx\times\Aa\mapsto S(\Bb)$ is  $\Aa$-lower semi-continuous.
\end{itemize}
\end{lemma}

To state the continuity theorems for minimax, 
we introduce
the multifunction  $\Phii_\Bb^{\Aa\leftrightarrow \Bb}:\Xx\times\Bb\mapsto 2^\Aa$ uniquely defined by its graph,
\begin{equation}\label{eq:auxil23}
\Gr(\Phii_\Bb^{\Aa\leftrightarrow \Bb}):=\{(x,b,a)\in \Xx\times\Bb\times\Aa\,:\, (x,a,b)\in \Gr(\Phii_{\Bb})\},
\end{equation}
that is, $\Phii_\Bb^{\Aa\leftrightarrow \Bb}(x,b)=\{a\in \Phii_\Aa(x)\,:\, b\in \Phii_\Bb(x,a)\},$ $(x,b)\in {\rm Dom\,}\Phii_\Bb^{\Aa\leftrightarrow \Bb}.$
We also introduce
the function $\ff^{\Aa\leftrightarrow \Bb}:\Gr(\Phii_\Bb^{\Aa\leftrightarrow \Bb})\subset (\Xx\times\Bb)\times\Aa\mapsto\R,$
\begin{equation}\label{eq:auxil3}
\ff^{\Aa\leftrightarrow \Bb}(x,b,a):=\ff(x,a,b),\quad (x,a,b)\in \Gr(\Phii_{\Bb}).
\end{equation}
According to (\ref{eq:auxil23}), the following equalities hold:
\begin{equation}\label{eq:auxil4}
\begin{aligned}
{\rm Dom\,}\Phii_\Bb^{\Aa\leftrightarrow \Bb}={\rm proj}_{\Xx\times\Bb}\Gr&(\Phii_{\Bb})=\{(x,b)\in \Xx\times\Bb\,:\,\\
&(x,a,b)\in  \Gr(\Phii_{\Bb}){\rm \ for\ some\ }a\in \Aa\},
\end{aligned}
\end{equation}
where ${\rm proj}_{\Xx\times\Bb}\Gr(\Phii_{\Bb})$ is a projection  of $\Gr(\Phii_{\Bb})$ on
$\Xx\times\Bb.$

We would like to mention that certain continuity properties of $\ff^\sharp,$ $\vv^\sharp,$ and $\Phii_\Bb^*$
do not use $\Aa$-lower semi-continuity of $\Phii_\Bb.$ In particular, the following statements hold: {\it \begin{itemize}
\item[\rm{(A)}] if  $\Phii_{\Bb}:\Gr(\Phii_{\Aa})\subset\Xx\times\Aa\mapsto S(\Bb)$ is lower semi-continuous multifunction and $\ff:{\rm Gr}(\Phii_{\Bb})\subset \Xx \times \Aa\times\Bb \mapsto \overline{\mathbb{R}}$ is lower semi-continuous function, then  $\ff^\sharp:
\Gr(\Phii_{\Aa})\subset\Xx\times\Aa\mapsto \overline{\mathbb{R}}$ is lower semi-continuous; \cite[Theorem~4]{FKZ2017};
\item[\rm{(B)}] if  $\ff:{\rm Gr}(\Phii_{\Bb})\subset (\Xx \times \Aa)\times\Bb \mapsto {\mathbb{R}}$ is $\K$-sup-compact on $\Gr(\Phii_{\Bb})$ and $\Phii_{\Aa}:\Xx\mapsto S(\Aa)$ is lower semi-continuous, then $\vv^\sharp:
\Xx\mapsto \mathbb{R}\cup\{-\infty\}$ and $\ff^\sharp:
\Gr(\Phii_{\Aa})\subset\Xx\times\Aa\mapsto \mathbb{R}\cup\{-\infty\}$ are upper semi-continuous; \cite[Theorems~6 and 9]{FKZ2017};
\item[\rm{(C)}] if additionally to assumptions from (A) and (B), $\Phii_{\Bb}:\Gr(\Phii_{\Aa})\subset\Xx\times\Aa\mapsto S(\Bb)$ is lower semi-continuous and  $\ff:\Gr(\Phii_\Bb)\subset \Xx\times\Aa\times\Bb\mapsto\mathbb{R}$ is lower semi-continuous, then $\Phii_{\Bb}^*:\Gr(\Phii_{\Aa})\subset\Xx\times\Aa\mapsto \K(\Bb)$ is upper semi-continuous; \cite[Theorems~4, 6 and 12]{FKZ2017}.
\end{itemize} }

The following theorem presents continuity results for the worst-loss function, minimax function, and solution multifunction $\Phii_\Aa^*$ that assume $\Aa$-lower semi-continuity of  $\Phii_{\Bb}:\Gr(\Phii_{\Aa})\subset\Xx\times\Aa\mapsto S(\Bb).$ 

\begin{theorem}\label{th:united}
Let $\Phii_{\Bb}:\Gr(\Phii_{\Aa})\subset\Xx\times\Aa\mapsto S(\Bb)$ be an $\Aa$-lower semi-continuous multifunction and the function $\ff^{\Aa\leftrightarrow \Bb}:\Gr(\Phii_\Bb^{\Aa\leftrightarrow \Bb})\subset (\Xx\times\Bb)\times\Aa\mapsto\mathbb{R}$ be $\K$-inf-compact on $\Gr(\Phii_\Bb^{\Aa\leftrightarrow \Bb}).$ Then the following statements hold: \begin{itemize}
\item[\rm(i)] the worst-loss function $\ff^\sharp:
\Gr(\Phii_{\Aa})\subset\Xx\times\Aa\mapsto \overline{\mathbb{R}}$ is $\K$-inf-compact on $\Gr(\Phii_{\Aa});$ Feinberg et al.~\cite[Theorem~5]{FKZ2017};
\item[\rm(ii)] the minimax function $\vv^\sharp:
\Xx\mapsto \overline{\mathbb{R}}$ is lower semi-continuous;
Feinberg et al.~\cite[Theorem~8]{FKZ2017};
\item[\rm(iii)] if additionally $\vv^\sharp:
\Xx\mapsto \mathbb{R}\cup\{-\infty\}$ is upper semi-continuous function (in view of (ii), $\vv^\sharp:
\Xx\mapsto \overline{\mathbb{R}}$ is lower semi-continuous and sufficient conditions for its upper semi-continuity are provided in statement~(B)), then the infimum in (\ref{eq1starminimax}) can be
replaced with the minimum and the solution multifunction $\Phii_{\Aa}^*:\Xx\mapsto S(\Aa)$ is upper semi-continuous and compact-valued; Feinberg et al.~\cite[Theorem~11]{FKZ2017}.
\end{itemize}
\end{theorem}

\begin{remark}
{\rm
Feinberg at el. \cite[Theorems~7, 10, and 13]{FKZ2017} contains additional results on continuity
properties of $\ff^\sharp,$ $\vv^\sharp,$ $\Phii_\Aa^*,$ and $\Phii_\Bb^*,$ which are combinations of Statements~(A,B,C) and Theorem~\ref{th:united} from above. These results include \cite[Theorem~13]{FKZ2017} described below before Example~\ref{exa:Alsc}.
}
\end{remark}
\section{Example}\label{s3}
In this section we provide an example demonstrating that the assumption that the multifunction $\Phii_{\Bb}:\Gr(\Phii_{\Aa})\subset\Xx\times\Aa\mapsto S(\Bb)$ is
$\Aa$-lower semi-continuous cannot be relaxed in each statement of Theorem~\ref{th:united} to the assumption that this multifunction is lower semi-continuous.

In the following  Example~\ref{exa:Alsc},
$\Phii_{\Aa}:\Xx\mapsto S(\Aa)$ and $\Phii_{\Bb}:\Gr(\Phii_{\Aa})\subset\Xx\times\Aa\mapsto S(\Bb)$ are continuous multifunctions, $\ff^{\Aa\leftrightarrow \Bb}:\Gr(\Phii_\Bb^{\Aa\leftrightarrow \Bb})\subset (\Xx\times\Bb)\times\Aa\mapsto\mathbb{R}$ is $\K$-inf-compact on $\Gr(\Phii_\Bb^{\Aa\leftrightarrow \Bb}),$ and $\ff:{\rm Gr}(\Phii_{\Bb})\subset (\Xx \times \Aa)\times\Bb \mapsto {\mathbb{R}}$ is $\K$-sup-compact on $\Gr(\Phii_{\Bb}),$
that is, all the assumptions of Statements~(A,B,C) and Theorem~\ref{th:united} hold, but  $\Phii_{\Bb}:\Gr(\Phii_{\Aa})\subset\Xx\times\Aa\mapsto S(\Bb)$ is not
$\Aa$-lower semi-continuous. Then none of statements~(i)--(iii) of Theorem~\ref{th:united} hold, that is, the worst-loss function $\ff^\sharp:
\Gr(\Phii_{\Aa})\subset\Xx\times\Aa\mapsto \overline{\mathbb{R}}$ is not
$\K$-inf-compact on $\Gr(\Phii_{\Aa}),$ the minimax function $\vv^\sharp:
\Xx\mapsto \mathbb{R}$ is not upper semi-continuous, and the solution multifunction $\Phii_{\Aa}^*:\Xx\mapsto S(\Aa)$ is not upper semi-continuous.

We recall that all the assumptions of statements~(A,B,C) and Theorem~\ref{th:united} taken together imply that the function $\ff^\sharp:
\Gr(\Phii_{\Aa})\subset\Xx\times\Aa\mapsto {\mathbb{R}}$ is continuous and $\K$-inf-compact on $\Gr(\Phii_{\Aa}),$ the function $\vv^\sharp:
\Xx\mapsto {\mathbb{R}}$ is continuous, and the multifunctins functions $\Phii_{\Aa}^*:\Xx\mapsto \K(\Aa)$ and $\Phii_{\Bb}^*:\Gr(\Phii_{\Aa})\subset\Xx\times\Aa\mapsto \K(\Bb)$ are upper semi-continuous;
Feinberg at el. \cite[Theorem~13]{FKZ2017}.

\begin{example}\label{exa:Alsc}
{\rm Let $\Xx:=\mathbb{R}$, $\Aa:=\Bb:=\mathbb{R}_+:=[0,+\infty),$
$\Phii_\Aa(x):=\mathbb{R}_+,$
$\Phii_\Bb(x,a):=[\phi_\Bb(x,a),+\infty),$
where
\[
\phi_\Bb(x,a):=
\left\{
\begin{array}{ll}
0, & \mbox{if either }x\le 0 \mbox{ or } x>0, 0\le a< \frac{1}{2x};\\ {}
2(2x+1)a-2- \frac1x, & \mbox{if }x>0 \mbox{ and } \frac{1}{2x} \le a \le \frac1x;
\\ {}
2+\frac1x , & \mbox{if }x>0 \mbox{ and } a>\frac1x;
\end{array}
\right.
\]
and let
\[
\ff(x,a,b):=
\left\{
\begin{array}{ll}
1+a-b, & \mbox{if either }x\le 0 \mbox{ or } x>0 ,\, 0\le a< \frac{1}{2x}, b\ge \phi_\Bb(x,a) ;\\ {}
(2x+1)a-b, & \mbox{if }x>0 ,\, \frac{1}{2x} \le a \le \frac1x, \mbox{ and }b\ge\phi_\Bb(x,a);
\\ {}
2+a-b, & \mbox{if }x>0 ,\, a>\frac1x, \mbox{ and }b\ge\phi_\Bb(x,a);
\end{array}
\right.
\] for all $x\in\Xx,$ $a\in \Phii_\Aa(x),$ and $b\in \Phii_\Bb(x,a).$

It is obvious that $\Phii_\Aa$ and  $\Phii_\Bb$ are continuous multifunctions because the constant function $x=0$ and the function $\phi_\Bb$ are continuous.

The multifunction $\Phii_\Bb$  is not $\Aa$-lower semi-continuous.  Indeed,
let 
$x:=a:=b:=0.$ Then $a\in \Phii_\Aa(x)$ and $b\in \Phii_\Bb(x,a).$
Let $x_n:=\frac1n \searrow x$ as $n\to +\infty$ and $a_n:=n\in \Phii_\Aa(x_n)=\mathbb{R}_+$ for all $n=1,2,\ldots.$  Then
$(-1,1)\cap \Phii_\Bb(x_{n},a_{n})=
\emptyset$ for each $n=1,2,\ldots.$ Therefore, $b=0$ is not a limit point of any sequence $\{b_n\}_{n=1,2,\ldots}$ with $b_n\in\Phii_\Bb(x_{n},a_{n}) ,$ $n=1,2,\ldots,$ because $|b_n-b|\ge1$ for each $n=1,2,\ldots,$ 
that is, $\Phii_\Bb$  is not $\Aa$-lower semi-continuous.

 In view of Lemma~\ref{k-inf-compact-strict}, the function $\ff^{\Aa\leftrightarrow \Bb}:\Gr(\Phii_\Bb^{\Aa\leftrightarrow \Bb})\subset (\Xx\times\Bb)\times\Aa\mapsto\mathbb{R}$ is $\K$-inf-compact on $\Gr(\Phii_\Bb^{\Aa\leftrightarrow \Bb})$ and the function $\ff:{\rm Gr}(\Phii_{\Bb})\subset (\Xx \times \Aa)\times\Bb \mapsto {\mathbb{R}}$ is $\K$-sup-compact on $\Gr(\Phii_{\Bb}).$ Therefore, these functions are continuous.

For every pair $(x,a)\in\Xx\times\Aa,$  the optimal decision for the second player is $b=\phi_\Bb(x,a).$
Thus,
\begin{equation}\label{eq:123}
\ff^\sharp(x,a)=\left\{
\begin{array}{ll}
1+a, & \mbox{if either }x\le 0 \mbox{ or } x>0, 0\le a< \frac{1}{2x} ;\\ {}
(2x+1)(\frac1x -a), & \mbox{if }x>0 \mbox{ and } \frac{1}{2x} \le a \le \frac1x;
\\ {}
a-\frac1x, & \mbox{if }x>0 \mbox{ and } a>\frac1x.
\end{array}
\right.
 \end{equation}
This function is continuous, but it is not $\K$-inf-compact on $\Xx\times\Aa$  because $x_n:=\frac1n\searrow0,$ the sequence $\ff^\sharp(\frac1n,n)=0,$ $n=1,2,\ldots,$ is bounded above, and $a_n:=n\to+\infty$ as $n\to\infty.$ Thus, the conclusion~(i) of Theorem~\ref{th:united} does not hold. The function $\ff^\sharp$ was introduced
in Luque-V\'asquez and Hern\'andez-Lerma \cite{LVHL}. In particular, \begin{equation}\label{eq:1234}
\begin{array}{lll}
\vv^\sharp(x)=
\left\{
\begin{array}{ll}
1, & \mbox{if }x\le 0  ;\\ {}
0, & \mbox{if }x>0;
\end{array}
\right.
&\mbox{ and }&
\Phii^*_\Aa(x)=
\left\{
\begin{array}{ll}
\{0\}, & \mbox{if }x\le 0  ;\\ {}
\{\frac1x\}, & \mbox{if }x>0;
\end{array}
\right.
\end{array}
\end{equation}
$x\in\Xx.$ The conclusions~(ii) and (iii) of Theorem~\ref{th:united} do not hold because the function $\vv^\sharp$ is not lower semi-continuous at $x=0$ and the solution multifunction $\Phii_\Aa^*$ is not upper semi-continuous at $x=0.$

In this example,  all the assumptions of statements~(A,B,C) and Theorem~\ref{th:united} hold except one: the multifunction $\Phii_\Bb$  is not $\Aa$-lower semi-continuous, but it is lower semi-continuous.  If the multifunction $\Phii_\Bb$ were $\Aa$-lower semi-continuous, then
the function $\ff^\sharp:
\Gr(\Phii_{\Aa})\subset\Xx\times\Aa\mapsto {\mathbb{R}}$ would be  $\K$-inf-compact on $\Gr(\Phii_{\Aa}),$ the function $\vv^\sharp:
\Xx\mapsto {\mathbb{R}}$ would be continuous, and the multifunction $\Phii_{\Aa}^*:\Xx\mapsto \K(\Aa)$ would be upper semi-continuous; see
Feinberg at el. \cite[Theorem~13]{FKZ2017}, whose description is provided before the example.
}
\end{example}

{\bf Acknowledgement.} Research of the first author was partially supported by NSF grant CMMI-1636193.



\medskip\noindent
{\bf \large {References}}

\end{document}